\documentclass[a4paper,12pt]{article}
\usepackage{latexsym}
\usepackage{amsfonts}
\usepackage{xypic}
\newtheorem{theorem}{Theorem}[section]
\newtheorem{corollary}[theorem]{Corollary}
\newtheorem{definition}{Definition}
\newtheorem{example}[theorem]{Example}
\newtheorem{remark}[theorem]{Remark}
\newtheorem{lemma}[theorem]{Lemma}
\newtheorem{proposition}[theorem]{Proposition}
\newtheorem{question}[theorem]{Question}
\title{Unification approach to the separation axioms between
$T_0$ and completely Hausdorff\thanks{1991 Math.\ Subject
Classification --- Primary: 54A05, 54D10; Secondary: 54D30,
54H05. \protect\newline Key words and phrases --- $\theta$-closed,
$\delta$-closed, Urysohn, zero-open, $\lambda$-closed, weakly
Hausdorff, $T_0$, $T_1$, completely Hausdorff, kc-space,
kd-space, $T_{\kappa,\lambda}$-space, anti-compact.}}
\author{Francisco G. Arenas\footnote{Research supported
by DGICYT grant PB95-0737.}, Julian Dontchev\footnote{Research
supported partially by the Ella and Georg Ehrnrooth Foundation
at Merita Bank, Finland.} and Maria Luz
Puertas\footnote{Corresponding author.}}
\begin{document}
\baselineskip=20pt plus 1pt minus 1pt
\newcommand{\fxy}{$f \colon (X,\tau) \rightarrow (Y,\sigma)$}
\maketitle
\begin{abstract}
The aim of this paper is to introduce a new weak separation axiom
that generalizes the separation properties between $T_1$ and
completely Hausdorff. We call a topological space $(X,\tau)$ a
$T_{\kappa,\xi}$-space if every compact subset of $X$ with
cardinality $\leq \kappa$ is $\xi$-closed, where $\xi$ is a
general closure operator. We concentrate our attention mostly
on two new concepts: kd-spaces and $T_{\frac{1}{3}}$-spaces.
\end{abstract}

\section{Introduction}\label{s1}

The definitions of most (if not all) weak separation axioms are
deceptively simple. However, the structure and the properties of
those spaces are not always that easy to comprehend.

In this paper we try to unify the separation axioms between $T_0$
and completely Hausdorff by introducing the concept of
$T_{\kappa,\xi}$-spaces. We call a topological space $(X,\tau)$
a $T_{\kappa,\xi}$-space if every compact subset of $X$ with
cardinality $\leq \kappa$ is $\xi$-closed where $\xi$ is a given
closure operator. With different settings on $\kappa$ and $\xi$
we derive most of the well-known separation properties `in the
semi-closed interval $[T_0,T_3)$'. We are going to consider not
only Kuratowski closure operators but more general closure
operators, such as the $\lambda$-closure operator \cite{ADG1} for
example ($\xi$ is a general closure operator on $X$ if $\xi
\colon {\rm exp}X \rightarrow {\rm exp}X$ and $\xi \emptyset =
\emptyset$, $A \subseteq \xi A$, $A \subseteq B \Rightarrow \xi
A \subseteq \xi B$, $\xi\xi A = \xi A$).

A subset $A$ of a topological space $(X,\tau)$ is called {\em
$\lambda$-closed} \cite{ADG1} if it is the intersection of a
closed set and a $\Lambda$-set (a {\em $\Lambda$-set} is a set
that is the intersection of a family of open sets \cite{M1}).
Complements of $\lambda$-closed sets are called {\em
$\lambda$-open}. The family of all $\lambda$-open sets is a
topology on $X$ if and only if $X$ is (by definition) a
$\lambda$-space \cite{ADG1}. This topology is coarser than $\tau$
and is denoted by ${\tau}_{\lambda}$.

For a topological space $(X,\tau)$, the family of all regular
open sets forms a base for a new topology $\tau_{s}$, coarser
than $\tau$, which is often called the {\em semi-regularization}
of $\tau$. A point $x \in X$ is called a {\em $\delta$-cluster
point} of a subset $A$ of $X$ \cite{V1} if $A \cap U \not=
\emptyset$ for every regular open set $U$ containing $x$. The set
of all $\delta$-cluster points of $A$ is called the {\em
$\delta$-closure} of $A$ and is denoted by ${\rm Cl}_{\delta}
(A)$. The set $A$ is called {\em $\delta$-closed} \cite{V1} if
$A = {\rm Cl}_{\delta} (A)$. Complements of $\delta$-closed sets
are called {\em $\delta$-open} and the family of all
$\delta$-open sets is a topology on $X$, coarser than $\tau$,
denoted by ${\tau}_{\delta}$.

The first part of the following lemma is well-known, while the
second part unites \cite[Lemma 4]{MRV1} and \cite[Lemma
1.1]{DG1}. Recall first that a subset $A$ of a space $(X,\tau)$
is called {\em locally dense} \cite{CM1} (= preopen) if $A
\subseteq {\rm Int}({\rm Cl}(A))$. Note that every open and every
dense set is locally dense. The family of all regular open
subsets of a topological space $(X,\tau)$ will be denoted by
$RO(X,\tau)$.

\begin{lemma}\label{l1}
{\rm (i)} ${\tau}_s = {\tau}_{\delta}$.

{\rm (ii)} If $A$ is a locally dense subset of a topological
space $(X,\tau)$, then:

\ \ \ \ {\rm (a)} $RO(A,\tau|A) = \{ R \cap A \colon R \in
RO(X,\tau) \}$.

\ \ \ \ {\rm (b)} $(\tau|A)_{s} = \tau_{s}|A$. $\Box$
\end{lemma}

Back in 1968, Veli\v{c}ko \cite{V1} introduced the concept of
$\theta$-open sets. Recall that set $A$ is called {\em
$\theta$-open} \cite{V1} if every point of $A$ has an open
neighborhood whose closure is contained in $A$. The {\em
$\theta$-interior} \cite{V1} of a $A$ in $X$ is the union of all
$\theta$-open subsets of $A$ and is usually denoted by ${\rm
Int}_{\theta} (A)$. Complement of a $\theta$-open set is called
{\em $\theta$-closed}. It is equivalent to stipulate that ${\rm
Cl}_{\theta} (A) = \{ x \in X \colon {\rm Cl}(U) \cap A \not=
\emptyset, U \in \tau$ and $x \in U \}$ and a set $A$ is
$\theta$-closed if and only if $A = {\rm Cl}_{\theta} (A)$. All
$\theta$-open sets form a topology on $X$, coarser than $\tau$,
usually denoted by $\tau_{\theta}$. Note that a space $(X,\tau)$
is regular if and only if $\tau = \tau_{\theta}$. Note also that
the $\theta$-closure of a given set need not be a $\theta$-closed
set -- however it is always $\delta$-closed. Moreover, always
$\tau_{\theta} \subseteq \tau_{\delta} \subseteq \tau$.

Recall that a topological space $(X,\tau)$ is called:

{\rm (1)} a {\em unique sequence space} (= $US$-space) if the
limit point of every converging sequence is uniquely determined,

{\rm (2)} a {kc-space} \cite{W1} if every compact set of $X$ is
closed.

{\rm (3)} {\em weakly Hausdorff} \cite{So1} if its
semi-regularization is $T_1$, i.e.\ if each singleton is
$\delta$-closed.

{\rm (4}) {\em h$T^R_1$ space} \cite{DS1} if every subspace of
$(X,\tau)$ is weakly Hausdorff.

A subset $A$ of a topological space $(X,\tau)$ is called {\em
zero-open} if for each $x \in A$ there exists a zero-set $Z$ (in
$X$) and a {\em cozero-set} $C$ (in $X$) such that $x \in C
\subseteq Z \subseteq A$. Complements of zero-open sets are
called {\em zero-closed} and the family of all zero-open sets is
a topology on $X$, coarser than $\tau$, denoted by ${\tau}_{z}$.
Recall that a topological space $(X,\tau)$ is called {\em
completely Hausdorff} if each two different points have disjoint
cozero neighborhoods.

A set $A$ is called {\em Urysohn-open} if for each $x \in A$
there exist two open sets $U$ and $V$ in $X$ such that $x \in U
\subseteq {\rm cl} (U) \subseteq V \subseteq {\rm cl} (V)
\subseteq A$. Complements of Urysohn-open sets are called {\em
Urysohn-closed} and the family of all Urysohn-open sets
${\tau}_{U}$ is a topology on $X$, coarser than $\tau$ and finer
than the quasi-topology $\tau_q$. The {\em quasi-topology}
$\tau_q$ on $X$ is the topology having as base the clopen subsets
of $(X,\tau)$. \cite{DGR1}.

\section{$T_{\kappa,\xi}$-spaces}\label{s2}

\begin{definition}\label{d1}
{\em A topological space $(X,\tau)$ is called a
{\em $T_{\kappa,\xi}$-space} if every compact subset of $X$ with
cardinality $\leq \kappa$ is $\xi$-closed, where $\xi$ is a
general closure operator. When $\xi$ is the usual closure
operator, we use the notation $T_{\kappa,c}$-space.}
\end{definition}

The following theorem shows that most of the well-known
definitions of separation axioms placed between $T_1$ and $T_2$
can be derived from the definition of $T_{\kappa,\xi}$-spaces.

\begin{theorem}
Let $(X,\tau)$ be a topological space with $|X| = \kappa$. Then:

{\rm (i)} $X$ is completely Hausdorff if and only $X$ it is a
$T_{\kappa,z}$-space.

{\rm (ii)} $X$ is Urysohn if and only if $X$ is a
$T_{\kappa,U}$-space.

{\rm (iii)} $X$ is Hausdorff if and only if $X$ is a
$T_{\kappa,\theta}$-space.

{\rm (iv)} $X$ is a kc-space if and only if $X$ is a
$T_{\kappa,c}$-space.

{\rm (v)} $X$ is weakly Hausdorff if and only if $X$ is a
$T_{1,\delta}$-space.

{\rm (vi)} $X$ is $T_1$ if and only if $X$ is a $T_{1,c}$-space.

{\rm (vii)} $X$ is $T_0$ if and only if $X$ is a
$T_{1,\lambda}$-space.

\end{theorem}

{\em Proof.} (i) Since $\{x\}$ is compact in $X$, it is
zero-closed, so  $X \setminus \{x\}$ is zero-open. Hence, if $y$
is a point of $X$ distinct from $x$ we have that $y \in X
\setminus \{x\}$, that is, there exist a zero set $Z$ and a
cozero set $C$ in $X$ such that $y \in C \subseteq Z \subseteq
X \setminus \{x\} $.Then $C$ and $X \setminus Z$ are disjoint
cozero neighborhoods of $x$ and $y$ respectively; hence the space
is completely Hausdorff.

Let $A$ be a compact subset of the completely Hausdorff space
$X$. Given a point $x \not\in A$ and a point $y \in A$, there
exists two cozero sets $C_y$ and $D_y$ with $x \in C_y$, $y \in
D_y$ and $C_y \cap D_y = \emptyset$. Clearly, $\{ D_y: y \in A\}$
is an open (cozero, in fact) covering of the compact set $A$, so
there exist a finite subcovering $\{ D_{y_{i}} \colon i=1,
\ldots, n \}$. Set $D = \bigcup_{i=1}^{n} D_{y_{i}}$. Note that
$D$ is a cozero set containing $A$. Now set $C =
\bigcap_{i=1}^{n} C_{y_{i}}$; $C$ is a cozero set containing $x$
and $C \cap D = \emptyset$; so we have $x \in C \subseteq X
\setminus D \subseteq X \setminus A$. As the point $x$ was
arbitrary and $X \setminus D$ is a zero set, we have that $X
\setminus A$ is zero-open, so $A$ is zero-closed, as desired.

(ii) Since $\{x\}$ is compact in $X$, it is Urysohn-closed, so
$X\setminus \{x\}$ is Urysohn-open. Hence, if $y$ is a point of
$X$ distinct from $x$, we have that $y \in X \setminus \{x\}$,
that is, there exist two open sets $U$ and $V$ in $X$ such that
$y \in U \subseteq {\rm cl} (U) \subset V \subset {\rm Cl} (V)
\subseteq X \setminus \{x\}$. Then, $U$ and $X \setminus {\rm Cl}
(V)$ are open neighborhoods of $x$ and $y$ respectively with
disjoint closures; hence the space is Urysohn.

Let $A$ be a compact subset of the Urysohn space $X$. Given
a point $x \not\in A$ and a point $y \in A$, there exists open
sets $U_y$ and $V_y$ with $x \in U_y$, $y \in V_y$ and ${\rm Cl}
(U_y) \cap {\rm Cl} (V_y) = \emptyset$. $\{ V_y \colon y \in A
\}$ is an open covering of the compact set $A$, so there exist
a finite subcovering $\{ V_{y_{i}} \colon i=1, \ldots, n \}$.
Take $V= \bigcup_{i=1}^{n} V_{y_{i}}$; $V$ is an open set
containing $A$ and ${\rm Cl} (V) = \bigcup_{i=1}^{n} {\rm Cl}
(V_{y_{i}})$. Now set $U = \bigcap_{i=1}^{n} U_{y_{i}}$. Clearly,
$U$ is an open set containing $x$ and ${\rm Cl} (U) =
\bigcap_{i=1}^{n} {\rm Cl} (U_{y_{i}})$. Note that ${\rm Cl} (U)
\cap {\rm Cl} (V) = \emptyset$; so we have $x \in U \subseteq
{\rm Cl} (U) \subset X \setminus {\rm Cl}(V) \subseteq X
\setminus V \subseteq X \setminus A$. As the point $x$ was
arbitrary, we have that $X \setminus A$ is Urysohn-open, so $A$
is
Urysohn-closed, as desired.

(iii) By Theorem 4.3 from \cite{J80}, a space is Hausdorff if and
only if every compact set is $\theta$-closed.

(iv) By definition a space is kc if and only if every compact set
is closed.

(v) A space is weakly Hausdorff if and only if its
semi-regularization is $T_1$ i.e.\ if every singleton is
$\delta$-closed \cite{So1}.

(vi) is obvious.

(vii) This is Theorem 2.5 from \cite{ADG1}. $\Box$

Of particular interest is probably the class of
$T_{\kappa,\delta}$-spaces, i.e.\ the spaces in which compact
sets are $\delta$-closed, since they have not been considered in
the literature until now. In order to be consistent with the
definition of kc-space we will call this class of spaces
kd-spaces. The relations between the spaces mentioned above are
given in the following diagram:

$$
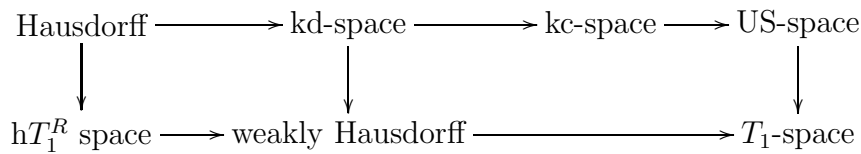
\diagram
\text{Hausdorff} \rto \dto & \text{kd-space} \rto \dto &
\text{kc-space} \rto & \text{US-space} \dto \\
\text{h$T^R_1$ space} \rto & \text{weakly Hausdorff} \rrto & &
\text{$T_1$-space}
\enddiagram
$$

\begin{example}\label{e1}
{\em Example of a weakly Hausdorff space which is not a kd-space,
not even a US-space: Let $\mathbb R$ be the real line and let $X
= [0,1] \cup [2,3] \cup S$, where $S = \{ 4,5,6, \ldots \}$. For
each $x \in X$ we define a neighborhood filter in the following
way. If $x \in [0,1) \cup (2,3] \cup S$, then the neighborhoods
of $x$ are the ones inherited by the usual topology on the real
line. A neighborhood of $1$ (resp.\ of $2$) are the sets
containing some interval $(a,1]$ (resp.\ $[2,b)$) along with all
but finitely many points of $S$. It is easily observed that we
have a topology on $X$. Since the sequence $(4,5,6, \ldots)$
converges both to $1$ and to $2$, then $X$ is not an $US$-space.
In particular $X$ is not a kd-space. Observe that both $1$ and
$2$ can be represented as intersection of regular closed sets,
hence they are both $\delta$-closed as is every other point of
$X$. Hence $X$ is weakly Hausdorff.}
\end{example}

\begin{example}\label{e2}
{\em Example of a kd-space which is not Hausdorff. Consider
Example
3 from \cite{SP1}. Let $X$ be the interval $[0,1]$ of real
numbers
with the following topology: all points from $(0,1)$ are clopen;
the basic neighbourhoods of $0$ are of the form $[0,x)$ for $x
>
0$; the neighbourhoods of $1$ are of the form $X \setminus F$,
where $F
\subseteq [0,1)$ is either a finite set or a sequence that
converges to $0$ with respect to the standard topology. It is
shown
in \cite{SP1} that this space is not Hausdorff. It can be easily
observed that all compact sets are closed, and hence
$\delta$-closed, since the space is semi-regular \cite{SP1}.}
\end{example}

Recall that a topological space $(X,\tau)$ is called a {\em
C'-space} if every compact set is compactly closed. A set $A
\subseteq X$ is called {\em compactly closed} if $A \cap K$ is
compact for any compact subset $K$ of $X$ or equivalently if the
canonical injection $i \colon A \rightarrow X$ is compact.

\begin{proposition}
Every kc and hence every kd-space is a C'-space.
\end{proposition}

\begin{example}
{\em The real line with the cofinite topology is an example of
a C'-space which is not a kc-space.}
\end{example}

In 1979, Bankston \cite{Ba1} introduced the anti operator on a
given topological space. Recall that a space $(X,\tau)$ is called
{\em anti-compact} if every compact subset of $X$ is finite.
Anti-compact spaces are sometimes called pseudo-finite. The class
of topological spaces where compact sets are finite was also
considered by in 1981 by Sharma \cite{Sh1}.

\begin{theorem}
For an anti-compact topological space $(X,\tau)$ the following
conditions are equivalent:

{\rm (1)} $X$ is a kd-space.

{\rm (2)} $X$ weakly Hausdorff.
\end{theorem}

{\em Proof.} Since (1) $\Rightarrow$ (2) is valid for any
topological space, then we only need to verify (2) $\Rightarrow$
(1). Let $K$ be a compact subset of $X$. Since $X$ is
anti-compact, then $K$ is finite. By (2), each point of $K$ is
$\delta$-closed and since the family of $\delta$-closed sets is
closed under finite additivity, then $K$ is $\delta$-closed.
Hence $X$ is a kd-space. $\Box$

\begin{example}\label{p2}
{\em There is a simple example of a kc-space which is not a
kd-space. Let $\mathbb R$ be the real line with the cocountable
topology $\tau_{cc}$. Since $({\mathbb R},\tau_{cc})$ is
anti-compact and $T_1$, it is clear that $({\mathbb
R},\tau_{cc})$
is a kc-space. But the semi-regularization topology is the
indiscrete one. Thus $({\mathbb R},\tau_{cc})$ is not a
kd-space.}
\end{example}

\begin{proposition}
{\rm (i)} For a first countable topological space $(X,\tau)$ the
following conditions are equivalent:

{\rm (1)} $X$ Hausdorff. \ {\rm (2)} $X$ is a kd-space. \
{\rm (3)} $X$ is a kc-space. \ {\rm (4)} $X$ is an US-space.

{\rm (ii)} A semi-regular topological space $(X,\tau)$ is a
kd-space if and only if it is a kc-space. $\Box$
\end{proposition}

\begin{remark}
{\em It is well-known that in the class of sequential spaces the
concepts of kc-spaces and US-spaces coincide. However, we do not
know if there exists an easy example of a sequential US-space
which is not a kd-space.}
\end{remark}

\begin{theorem}
Locally dense subspaces of kd-spaces are kd-spaces.
\end{theorem}

{\em Proof.} Follows easily from Lemma~\ref{l1} (ii). $\Box$

\begin{corollary}
Let $(X_{\alpha},{\tau}_{\alpha})_{\alpha \in \Omega}$ be a
family of topological spaces. For the topological sum $X =
\sum_{\alpha \in \Omega} X_{\alpha}$ the following conditions are
equivalent:

{\rm (1)} $X$ is a kd-space.

{\rm (2)} Each $X_{\alpha}$ is a kd-space. $\Box$
\end{corollary}

According to Mr\v{s}evi\'{c}, Reilly and Vamanamurthy a function
\fxy\ is called {\em super-closed} (resp.\ 
{\em super-continuous}) \cite{MRV1} if the image of every closed
subset of $X$ is $\delta$-closed in $Y$ (resp.\ the preimage of
every open set is $\delta$-open). A bijection \fxy\ is called
{\em a super-homeomorphism} if $f$ is both super-closed and
super-continuous.

\begin{proposition}
If $(X,\tau)$ is compact and $(Y,\sigma)$ is a kd-space, then
every continuous function \fxy\ is perfect and super-closed and
hence every (super-)continuous bijection \fxy\ is a
(super-)homeomorphism.
\end{proposition}

{\em Proof.} Let $K$ be a closed subset of $X$. Clearly $K$ is
compact. Since $f$ is continuous, then $f(K)$ is compact. Since
$Y$ is a kd-space, then $f(K)$ is $\delta$-closed and hence $f$
is super-closed. For any $y \in Y$, $\{ y \}$ is $\delta$-closed,
since $Y$ is a kd-space. Since $f$ is continuous, then $f^{-1}
\{ y \}$ is closed and moreover compact, since $X$ is compact.
Thus $f$ is perfect. The rest of the claim is obvious. $\Box$

\section{$T_{\kappa,\lambda}$-spaces}\label{s3}

Yet another class of particular interest is probably the class
of {\em $T_{\kappa,\lambda}$-spaces}, i.e.\ the spaces in which
compact sets are $\lambda$-closed, as they have not been also
considered in the literature until now. In order to be consistent
with the fact that they are placed between $T_{\frac{1}{2}}$
(spaces where every set is $\lambda$-closed, see Theorem 2.6 of
\cite{ADG1}) and $T_{\frac{1}{4}}$ (spaces where every finite set
is $\lambda$-closed, see Theorem 3.1 of \cite{ADG1}), we call
them $T_{\frac{1}{3}}$ spaces. Are they properly placed between
the other two separation axioms? The following examples shows
that they are.

\begin{example}\label{ea1}
{\em Let $X$ be the set of non-negative integers with the
topology whose open sets are those which contain $0$ and have
finite complement (so closed sets are the finite sets that do not
contain $0$). Every point is closed except $0$, which is neither
open nor closed nor even locally closed. This space is neither
$T_{\frac{1}{2}}$ nor $T_{D}$ (= singletons are locally closed),
although it is $T_0$. However it is $T_{{\frac{1}{4}}}$. To see
that it is not $T_{\frac{1}{3}}$, just note that every subset of
$X$ is compact, so if it is $T_{\frac{1}{3}}$, it is also
$T_{\frac{1}{2}}$, but we have proved it is not.}
\end{example}

\begin{example}\label{ea2}
{\em Next we present an example of a $T_{\frac{1}{3}}$-space that
is not $T_{\frac{1}{2}}$. Let $X$ be an uncountable set and let
$p$ a fixed point in $X$. Then the family $\tau_{p} = \{
\emptyset \} \cup \{ U \subset X : p \in U$ and $X \setminus U$
is countable$\}$ is called the cocountable topology in $X$
generated by $p$ and any space $X$ equipped with such a topology
is called {\em cocountably point generated space}. Using
Theorem~\ref{ttfin1}, we first show that this space is
$T_{\frac{1}{3}}$. Let $K$ be compact and $y \not\in K$. If $y
\not= p$, then we are done, since $\{ y \}$ is closed. If $y =
p$, then $K$ does not contain $y$. We will show that $K$ is
finite. If $K$ is infinite, let $K = A \cup B$, where $A$ and $B$
are disjoint and $A$ is denumerable. Now, $\{ B \cup \{ x \}
\colon x \in A \}$ is an infinite open cover of $(K,\tau_{p}|K)$,
which has no finite subcover. By contradiction, $K$ is finite and
thus closed. However, $(X,\tau_{p})$ is not a
$T_{\frac{1}{2}}$-space, since $\{ p \}$ is neither open nor
closed.}
\end{example}

The proof of the following result is easy and hence omitted.

\begin{proposition}
For an anti-compact topological space $(X,\tau)$ the following
conditions are equivalent:

{\rm (1)} $X$ is $T_{\frac{1}{3}}$.

{\rm (2)} $X$ is $T_{\frac{1}{4}}$.
\end{proposition}

\begin{question}
{\em The only compact subsets of the space in Example~\ref{ea2}
are the finite ones. Is it true that if a topological space $X$
is $T_{\frac{1}{3}}$ and is not $T_{\frac{1}{2}}$, then $X$ is
anti-compact?}
\end{question}

We have the following characterization of
$T_{\frac{1}{3}}$-spaces.

\begin{theorem}\label{ttfin1}
For a topological space $(X,\tau)$ the following conditions are
equivalent:

{\rm (1)} For every compact subset $F$ of $X$ and every $y
\not\in F$ there exists a set $A_y$ containing $F$ and disjoint
from $\{ y\}$ such that $A_y$ is either open or closed.

{\rm (2)} $X$ is $T_{\frac{1}{3}}$.
\end{theorem}

{\em Proof.} (1) $\Rightarrow$ (2) Let $F \subseteq X$ be a
compact subset of $X$. Then for every point $y \not\in F$ there
exists a set $A_y$ containing $F$ and disjoint from $\{ y \}$
such that $A_y$ is either open or closed. Let $L$ be the
intersection of all open sets $A_y$ and let $C$ be the
intersection of all closed sets $A_y$. Clearly, $L$ is a
$\Lambda$-set and $C$ is closed. Note that $F = L \cap C$. This
shows that $F$ is $\lambda$-closed, and hence is
$T_{\frac{1}{3}}$.

(2) $\Rightarrow$ (1) Let $F$ be a compact subset of $X$ and $y$
be a point of $X \setminus F$. Since $X$ is $T_{\frac{1}{3}}$ $F$
is $\lambda$-closed, $F = L \cap C$, where $L$ is a $\Lambda$-set
and $C$ is closed. If $C$ does not contain $y$, then $X \setminus
C$ is an open set containing $y$ and we are done. If $C$ contains
$y$, then $y \not\in L$ and thus for some open set $U$,
containing $F$, we have $y \not\in U$, and we are also done.
$\Box$

Recall that topological space $(X,\tau)$ is called a {\em weak
$R_0$-space} \cite{DGR1} if every $\lambda$-closed singleton is
a $\Lambda$-set. Note that $T_{\frac{1}{3}}$ neither implies weak
$R_0$ nor $R_0$, since a space is $T_1$ if and only if is $T_0$
and $R_0$ if and only if is $T_0$ and weak $R_0$ (Theorem 2.9 of
\cite{ADG1}), and we have that $T_{\frac{1}{3}}$ implies $T_0$
and there are $T_{\frac{1}{3}}$ spaces that are not $T_1$.

\begin{question}
Is there a (nice) characterization for semiregular spaces in
terms of $T_{\kappa,\xi}$-spaces?
\end{question}

{\bf Acknowledgement.} The authors thank the referee for his help
in improving the quality of this paper.

\
\
\begin{center}
Area of Geometry and Topology\\Faculty of
Science\\Universidad de Almer\'{\i}a\\04071
Almer\'{\i}a\\Spain\\e-mail: {\tt farenas@ualm.es}
\end{center}
\
\begin{center}
Department of Mathematics\\University of Helsinki\\PL 4,
Yliopistonkatu 15\\00014 Helsinki\\Finland\\e-mail: {\tt
dontchev@cc.helsinki.fi}, {\tt
dontchev@e-math.ams.org}\\http://www.helsinki.fi/\~{}dontchev/
\end{center}
\
\begin{center}
Area of Geometry and Topology\\Faculty of
Science\\Universidad de Almer\'{\i}a\\04071
Almer\'{\i}a\\Spain\\e-mail: {\tt mpuertas@ualm.es}
\end{center}
\
\
\end{document}